\documentclass{amsart}
\usepackage{amssymb}
\usepackage{a4wide}
\usepackage[english]{babel}

\newcommand{\Ps}{\mathbb{P}}
\newcommand{\Z}{\mathbb{Z}}
\newcommand{\C}{\mathbb{C}}
\newcommand{\Li}{\mathbb{L}}
\newcommand{\ra}{\rightarrow}

\renewcommand{\phi}{\varphi}
    \newtheorem{Lem}{Lemma}[section]
    \newtheorem{Prop}[Lem]{Proposition}
    \newtheorem{Thm}[Lem]{Theorem}
    
    \newtheorem{Cor}[Lem]{Corollary}
    \newtheorem{Que}[Lem]{Question}
    \newtheorem{Ass}[Lem]{Assumption}

   \theoremstyle{definition}
    \newtheorem{Def}[Lem]{Definition}
    \newtheorem{Exa}[Lem]{Example}
    \newtheorem{Rem}[Lem]{Remark}
    
    \DeclareMathOperator{\rank}{rank}
    \DeclareMathOperator{\Hom}{Hom}

\begin{document}
\title{Extremal elliptic surfaces \& Infinitesimal Torelli}
\author{Remke Kloosterman}

\address{Department of Mathematics and Computer Science, University of Groningen, PO Box 800, 9700 AV  Groningen, The Netherlands}
\email{r.n.kloosterman@math.rug.nl}
\thanks{Part of this research was done during the author's stay as EAGER pre-doc at the Turin node of EAGER. The author would like to thank Alberto Conte and Marina Marchisio for making this possible. The author wishes to thank Bert van Geemen and Jaap Top for the useful conversations on this topic and to thank  
Frederic Mangolte for suggesting \cite{Nori}. The author wishes to thank the referee for suggesting several improvements. Parts of this paper will be part of the author's PhD thesis.}
\begin{abstract}
  We describe in terms of the $j$-invariant all elliptic surfaces $\pi: X \rightarrow C$ with a section, such that $h^{1,1}(X)=\rank NS(X)$ and the Mordell-Weil group of $\pi$ is finite.

We use this to give a complete solution to infinitesimal Torelli for elliptic surfaces over $\Ps^1$ with a section.
\end{abstract}
\subjclass{ 14J27  (Primary); 14J28, 14D07, 32G20 (Secondary)}

\keywords{Elliptic surfaces, Infinitesimal Torelli}
\date{\today}
\maketitle

\section{Introduction}\label{Intro}
An extremal elliptic surface over $\C$ is an elliptic surface such that the rank of the N\'eron-Severi group equals $h^{1,1}$ and there are finitely many sections. 

Their main application, until now, is in the classification of singular fibers on certain elliptic surfaces: Once a configuration of singular fibers on an extremal elliptic surface is known, one can construct from this configuration many other configurations of singular fibers on elliptic surfaces, where the genus of the base curve and the geometric and arithmetic genus of the surface remain fixed. In \cite{Sh} these operations are called elementary transformations, and are, a priori, only valid for $K3$ surfaces. Actually, all elementary transformations are combinations of twisting and deformations of the $J$-map (terminology from \cite[VIII.2]{Mi1} and \cite{Mi2}) and hence are valid for any elliptic surface.

The classification of singular fibers on a rational elliptic surface has been given more than 10 years ago (see \cite{Per},  \cite{Mi2}, \cite{OS}). Recently there has been given a classification of all singular fibers of elliptic $K3$ surfaces with a section (see \cite{Sh}). From the classification of configurations of singular fibers on rational (see \cite{Mi2}) and $K3$ surfaces (see \cite{Sh}) we know that any configuration can be obtained from an extremal configuration using elementary transformations.

In this paper we give a complete classification of extremal elliptic surfaces with constant $j$-invariant (Theorem~\ref{constthm}). We use this classification to prove the following.
\begin{Thm}\label{torthm} Let $\pi: X \ra \Ps^1$ be an elliptic surface without multiple fibers. Assume that $p_g(X)>1$. Then $X$ does not satisfy infinitesimal Torelli if and only if $j(\pi)$ is constant and $\pi$ is extremal. \end{Thm}

Ki\u\i \; (\cite[Theorem 2]{Kii}) proved infinitesimal Torelli for elliptic surfaces without multiple fibers and non-constant $j$-invariant. Sait\=o (\cite{Sai}) proved in a different way infinitesimal Torelli for elliptic surfaces without multiple fibers and $j$-invariant different from 0 and 1728.

For elliptic surfaces with non-constant $j$-invariant we will give the following structure theorem:

\begin{Thm} \label{mthm}Suppose $\pi:X\rightarrow C$ is an elliptic surface without multiple fibers and non-constant $j$-invariant, then the following are equivalent:
\begin{enumerate}
\item $\pi$ is extremal 
\item $j(\pi)$ is not ramified outside $0, 1728, \infty$, the only possible ramification indices above $0$ are $1,2,3$ and above $1728$ are $1,2$, and $\pi$ has no fibers of type $II, III, IV$ or $I_0^*$.
\item There exists an elliptic surface with $\pi': X' \ra C$, with $j(\pi')=j(\pi)$, the fibration $\pi'$ has no fibers of type $II^*,III^*$ or $IV^*$, at most one fiber of type $I_0^*$, and $\pi'$ has precisely $2p_g(X)+4-4g(C)$ singular fibers.
\end{enumerate}
\end{Thm}

We will give a similar theorem which includes the Mordell-Weil group. Let $m,n \in \Z_{\geq 1}$ be such that $m|n$ and $n>1$. Let $X_m(n)$ be the modular curve parameterizing triples $((E,O),P,Q)$, such that $(E,O)$ is an elliptic curve, $P\in E$ is a point of order $m$ and $Q\in E$ a point of order $n$.

If $(m,n)\not \in \{ (1,2), (2,2), (1,3), (1,4), (2,4) \}$ then there exists a universal family for $X_m(n)$, which we denote by $E_m(n)$. Denote by $j_{m,n}: X_m(n) \rightarrow \Ps^1$ the map usually called $j$.

From the results of \cite[\S 4 \& \S 5]{ShEMS} it follows that $E_m(n)$ is an extremal elliptic surface. The following theorem explains how to construct many examples of extremal elliptic surfaces with a given torsion group.

\begin{Thm}\label{mwthm} Fix $m,n\in \Z_{\geq 1}$ such that $m|n$ and $(m,n)\not \in \{(1,1), (1,2), (2,2), (1,3), (1,4),(2,4) \}$. Let $j\in \C(C)$ be non-constant then there exists a unique elliptic surface $\pi:X \rightarrow C$, with $j(\pi)=j$ and $MW(\pi)$ has $\Z/n\Z\times \Z/m\Z$ as a subgroup if and only if $j$ is of $(3,2)$-type, not ramified outside $0,1728, \infty$ and $j=j_{n,m} \circ g$ for some $g: C \rightarrow X_m(n)$.

The extremal elliptic surface with $j(\pi)=j$ and $\Z/n\Z \times \Z/m\Z$ is a subgroup of $MW(\pi)$ is the unique surface with only singular fibers of type $I_\nu$.
\end{Thm}

If $\pi: X \ra \Ps^1$ is an extremal semi-stable rational elliptic surface  then $X$ is determined by the configuration of singular fibers (see \cite[Theorem 5.4]{MP}). It seems that this quite particular for rational elliptic surfaces. If $X$ is a $K3$ surface a similar statement does not hold:

\begin{Thm}\label{unthm} There exists pairs of {\em extremal semi-stable elliptic $K3$ surfaces}, $\pi_i: X_i \rightarrow \Ps^1$, $(i=1,2)$ such that $MW(\pi_1)\cong MW(\pi_2)$, the configuration of singular fibers of the $\pi_i$ coincide and $X_1$ and $X_2$ are non-isomorphic. \end{Thm}

This gives a negative answer to \cite[Question 0.2]{ATZ}. The essential ingredient for the proof comes from \cite[Table 2]{SZ}.

The paper is organized as follows:

Section~\ref{Prel} contains some definitions and several standard facts. In section~\ref{cst} we give a list of extremal elliptic surfaces with constant $j$-invariant. They behave different from the non-constant ones. There are exactly 5 infinite families of extremal elliptic surfaces with constant $j$-invariant (3 of dimension 1, 1 of dimension 2 and 1 of dimension 3). In section~\ref{Tor} we explain this different behavior by proving Theorem~\ref{torthm}. In section~\ref{twist} we explain how twisting can reduce the problem of classification. In section~\ref{singfib} we link the ramification of the $j$-map and the number of singular fibers of a certain elliptic surface. This combined with the results of section~\ref{twist} gives a proof of Theorem~\ref{mthm}. In section~\ref{MW} contains a proof of the version with the description of the group of sections (Theorem~\ref{mwthm}). In section~\ref{uni} we prove Theorem~\ref{unthm}. In section~\ref{cla} we give a classification of extremal elliptic surfaces with $g(C)=p_g(X)=q(X)=1$. Section~\ref{one} contains a proof of the fact that for any positive $k$ there exist elliptic surfaces with only one singular fiber; the singular fiber has to be of type $I_{12k}$ or $I_{12k-6}^*$, and both occur.
\section{Preliminaries and conventions}
\label{Prel}
\begin{Ass}By a curve we mean a non-singular projective complex curve.

By a surface we mean a non-singular projective complex surface.\end{Ass}

\begin{Def}\label{defbas} An \emph{elliptic surface} is a triple $(\pi,X,C)$ with $X$ a surface, $C$ a curve, $\pi$ is a morphism  $X\rightarrow C$, such that almost all fibers are irreducible genus 1 curves and $X$ is relatively minimal.

  We denote by $j(\pi): C \rightarrow \Ps^1$ the rational function such that $j(\pi)(P)$ equals the $j$-invariant of $\pi^{-1}(P)$, whenever $\pi^{-1}(P)$ is non-singular.

A \emph{Jacobian elliptic surface} is an elliptic surface together with a section $\sigma_0: C \rightarrow X$ to $\pi$.

  The set of sections of $\pi$ is an abelian group, with $\sigma_0$ as the identity element. Denote this group by $MW(\pi)$.

  An \emph{extremal elliptic surface} is an elliptic surface such that $\rho(X)=h^{1,1}(X)$ and $MW(\pi)$ is empty or finite.

  Let $\Li$ be the fundamental line bundle $[R^1 \pi_* \mathcal{O}_X]^{-1}$.
\end{Def}

\begin{Ass}\label{mainass} All elliptic surface are without multiple fibers.\end{Ass}

\begin{Rem} To an elliptic surface $\pi:X \ra C$ we can associate its Jacobian fibration $Jac(\pi): Jac(X) \ra C$. The Hodge numbers $h^{p,q}$, the Picard number $\rho(X)$, the type of singular fibers of $\pi$ and $\deg(\Li)$ are invariant under taking the Jacobian fibration.\end{Rem}

\begin{Rem} If $P$ is a point on $C$, such that $\pi^{-1}(P)$ is singular then $j(\pi)(P)$ behaves as follows:
\[ \begin{array}{|c|c|}
\hline
\mbox{Kodaira type of fiber over }P & j(\pi)(P) \\
\hline
I_0^*                  & \neq \infty \\
I_\nu, I_\nu^* (\nu>0) & \infty \\
II, IV, IV^*, II^*     & 0  \\
III,III^*              & 1728\\ \hline \end{array}\]
\end{Rem}

\begin{Def} Let $X$ be a surface, let $C$ and $C_1$ be curves. Let $\varphi: X \rightarrow C$ and $f: C_1 \rightarrow C$ be two morphisms. Then we denote by $X \times_C C_1$ the smooth, relatively minimal model of the ordinary fiber product of $X$ and $C_1$.
\end{Def}

We use the line bundle $\Li$ only to simplify notation. Note that $\deg(\Li)=p_g(X)+1-g(C)=p_a(X)+1$. (See \cite[Lemma IV.1.1]{Mi1}.)

Recall the following theorem.
\begin{Thm}[{Shioda (\cite[Theorem 1.3 \& Corollary 5.3]{Sd})}]\label{ST}  Let $\pi:X\rightarrow C$ be a Jacobian elliptic surface, such that $\deg(\Li)>0$. Then the N\'eron-Severi group of $X$ is generated by the classes of $\sigma_0(C)$, a fiber, the components of the singular fibers not intersecting $\sigma_0(C)$, and the generators of the Mordell-Weil group. Moreover, let $S$ be the set of points such that $\pi^{-1}(P)$ is singular. Let $m(P)$ be the number of irreducible components of $\pi^{-1}(P)$, then
  \[ \rho(X) := \rank(NS(X))=2+ \sum_{P \in S} (m(P)-1)+\rank(MW(\pi)) \]
\end{Thm}

\begin{Def} Suppose $\pi: X \ra C$ is an elliptic fibration. Denote by $\Lambda(Jac(\pi))$ the subgroup of the N\'eron-Severi group of $Jac(\pi)$ generated by the classes of the fiber, $\sigma_0(C)$ and the components of the singular fibers not intersecting $\sigma_0(C)$. Let $\rho_{tr}(\pi)=\rank \Lambda(Jac(\pi))$. 
\end{Def}

Note that if $\pi:X \ra \Ps^1$ has $\deg(\Li)=0$, then $\rho_{tr}=2$, although Theorem~\ref{ST} does not apply.

\begin{Def} Let $\pi: X \rightarrow C$ be an elliptic surface, define
  \begin{itemize}
\item $a(\pi)$ as the number of fibers of type $II^*,III^*,IV^*$.
\item $b(\pi)$ as the number of fibers of type $II,III,IV$.
\item $c(\pi)$ as the number of fibers of type $I_0^*$.
\item $d(\pi)$ as the number of fibers of type $I_\nu^*$, with $\nu>0$.
\item $e(\pi)$ as the number of fibers of type $I_\nu$, $\nu>0$.
\end{itemize}
\end{Def}

\begin{Def} Let $\pi: X \rightarrow C$ be an elliptic surface. Let $P\in C(\C)$. Define $v_P(\Delta_P)$ as the valuation at $P$ of the minimal discriminant of the Kodaira-N\'eron model, which equals the topological Euler characteristic of $\pi^{-1}(P)$. \end{Def}

\begin{Prop}\label{Noether} Let $\pi: X \rightarrow C$ be an elliptic surface. Then
\[ \sum_{P\in C(\C)} v_P(\Delta_P) = 12\deg(\Li) \]
\end{Prop}

\begin{proof} This follows from Noether's formula (see \cite[III.4.4]{Mi1}).
\end{proof}

\begin{Prop} \label{fibcrit} For any elliptic surface $\pi:X \rightarrow C$, not a product, we have
  \[h^{1,1}(X)-\rho_{tr}(\pi)=2(a(\pi)+b(\pi)+c(\pi)+d(\pi))+e(\pi)-2\deg(\Li)-2+2g(C)\]
\end{Prop}

\begin{proof} Recall from \cite[Lemma IV.1.1]{Mi1} that
\[h^{1,1}=10\deg(\Li)+2g(C).\]
From Kodaira's classification of singular fibers and Proposition~\ref{Noether} it follows that
\[ \rho_{tr}(\pi)=2+12\deg(\Li)-2(a(\pi)+b(\pi)+c(\pi)+d(\pi))-e(\pi).\]
Combining these finishes the proof.
\end{proof}

\begin{Cor} \label{cstcrt} Let $\pi:X \rightarrow C$ be an elliptic surface with constant $j$-invariant, not a product. Then $\pi$ is extremal if and only if  $\pi$ has $\deg(\Li)+1-g(C)$ singular fibers.
\end{Cor}
\begin{proof} If $j$ is constant then $e(\pi)=d(\pi)=0$.\end{proof}

\section{Constant $j$-invariant}
\label{cst}
In this section we give a list of all extremal elliptic surfaces with constant $j$-invariant. 

\begin{Lem}~\label{curvlem}
Suppose $\pi: X \rightarrow C$ is an extremal elliptic surface such that $j(\pi)$ is constant. Then $g(C)\leq 1$.\end{Lem}

\begin{proof} If $j(\pi)$ is constant $v_P(\Delta_P)\leq 10$, for every point $P$. From this it follows from Proposition~\ref{Noether} and Corollary~\ref{cstcrt} that
\[ 12 \deg(\Li) = \sum_{P | \pi^{-1}(P) \mbox{ singular}} v_P(\Delta_p) \leq 10 (\deg(\Li)+1-g(C)). \]
From this it follows that $g(C)\leq 1$.
\end{proof}

\begin{Lem}~\label{classlem}
Suppose $\pi: X \rightarrow C$ is an extremal elliptic surface such that $j(\pi)$ is constant. Then one of the following occurs
\begin{enumerate}
\item $g(C)=1$; $\deg(\Li)=0$ and $X$ is a hyperelliptic surface.
\item $g(C)=0$; $j(\pi) \neq 0, 1728$;  $\deg(\Li)=1$.
\item $g(C)=0$; $j(\pi)=0$;  $1 \leq \deg(\Li)\leq 5$.
\item $g(C)=0$; $j(\pi)=1728$;  $1 \leq \deg(\Li) \leq 3$.
\end{enumerate}
\end{Lem}

\begin{proof} Suppose $g(C)=1$. Then $\pi$ has $\deg(\Li)$ singular fibers so
\[ 12 \deg(\Li) = \sum v_P(\Delta_P) \leq 10 \deg(\Li) \]
from which it follows that $\deg(\Li)=0$, so $X$ is a hyperelliptic surface.

Suppose $g(C)=0$. If $\deg(\Li)=0$ then $\pi$ is birational to a projection from a product, so by definition $\pi$ is not extremal.

Suppose $\deg(\Li)>0$. Assume $j(\pi)\neq 0, 1728$. Then all singular fibers are of type $I_0^*$. Since the Euler characteristic of such a fiber is 6, Proposition~\ref{Noether} implies that there are exactly $2\deg(\Li)$ singular fibers. Applying Corollary~\ref{cstcrt} gives
  \[ 2\deg(\Li)=\deg(\Li)+1 \]
  From this we know $\deg(\Li)=1$.

  Assume $j(\pi)=1728$ then all singular fibers are of type $III, I_0^*, III^*$. From this it follows that $v_P(\Delta_P)\leq 9$. By Proposition~\ref{Noether} and Proposition~\ref{fibcrit}
  \[12\deg(\Li) \leq 9 (a(\pi)+b(\pi)+c(\pi)) = 9\deg(\Li)+9\]
  so $1 \leq \deg(\Li) \leq 3$.

  Assume $j(\pi)=0$ then similarly we obtain $1\leq \deg(\Li) \leq 5$.
\end{proof}

\begin{Thm}\label{constthm} Suppose $\pi: X \rightarrow C$ is an elliptic surface with $j(\pi)$ constant.

Then $\pi$ is extremal if and only if either $C$ is a curve of genus 1 and $Jac(X)$ is a hyperelliptic surface or $C\cong \Ps^1$ and $Jac(\pi)$ has a model isomorphic to one of the following:

\begin{itemize}
\item $(j(\pi)=0) \; y^2=x^3+f(t)$ where $f(t)$ comes from the following table (the left hand side indicates the positions of the singular fibers)
\[\begin{array}{|c|c|c|c|c||c|l|}
\hline
  II &IV & I_0^*&  IV^*& II^*&p_g&f(t)\\
\hline
0 &&&&              \infty&0 &t\\
 &  0  & & \infty& & 0 &t^2\\
 & & 0,\infty & & & 0& t^3\\

 &  1   &&&             0,\infty& 1 & t^5(t-1)^2\\
 & &       1    &  0  &    \infty&1& t^4(t-1)^3\\
 &&& 0,1,\infty & & 1&t^4(t-1)^4\\

 && \alpha&&0,1,\infty & 2& t^5(t-1)^5 (t-\alpha)^3\\
 &&&    \alpha,1&  0,\infty & 2 & t^5(t-1)^4 (t-\alpha)^4\\

 &&&\beta&\alpha, 0,1,\infty  & 3 &  t^5(t-1)^5 (t-\alpha)^5(t-\beta)^4\\

 &&&&0,1,\infty,\alpha,\beta,\gamma & 4&t^5(t-1)^5 (t-\alpha)^5(t-\beta)^5(t-\gamma)^5 \\
 \hline
\end{array}\]
where $\alpha,\beta,\gamma \in \C-\{0,1\}$, pairwise distinct.

\item $(j(\pi)=1728) \; y^2=x^3+g(t)x$ where $g(t)$ comes from the following table
\[\begin{array}{|c|c|c||c|l|}
  \hline
  III&I_0^*& III^*& p_g &g(t)\\
  \hline
  0&  & \infty& 0 & t\\
  & 0,\infty  & &0&t^2 \\
  &1&  0,\infty & 1&t^3(t-1)^2\\
  & & 0,1,\infty, \alpha & 2&t^3 (t-1)^3(t-\alpha)^3 \\
  \hline
  \end{array}\]
where $\alpha \neq 0,1$.

\item $(j(\pi)\neq 1728)\; y^2=x^3 + a t^2 x + t^3$, with singular fibers of type $I_0^*$ at $t=0$ and $t=\infty$
\end{itemize}
\end{Thm}

\begin{proof} The above list follows directly from Corollary~\ref{cstcrt} and Lemma~\ref{classlem}. Since all cases are very similar, we discuss only the case $j(\pi)=0$ and $p_g=2$. In this case $Jac(\pi)$ has a model isomorphic to
\[ y^2=x^3+f(t) \]
with $f$ a polynomial such that $13 \leq \deg(f) \leq 18$, and $v_P(f)\leq 5$ for all finite $P$. At all zeros of $f$ there is a singular fiber. If $\deg(f)<18$ then the fiber over $t=\infty$ is also singular. 

If $\pi$ is extremal then from Corollary~\ref{cstcrt} it follows that $\pi$ has exactly 4 singular fibers. Assume that the fibers with the highest Euler characteristic are over $t=\infty, 0 ,1$.  Since $5+5+5+3=5+5+4+4$ are the only two ways of writing 18 as a sum of four positive integers smaller then 6, we obtain that after applying an isomorphism, if necessary, $f$ equals either 
\[t^5(t-1)^5(t-\alpha)^3 \mbox{ or }  t^5(t-1)^4(t-\alpha)^4. \]
\end{proof}

\begin{Rem}\label{modrem} Note that all extremal elliptic surfaces with constant $j$-invariant and $p_g(X)>1$ have moduli.\end{Rem}

 \section{Infinitesimal Torelli} \label{Tor}
In the previous section we gave examples of families of elliptic surfaces with maximal picard number. In this section we prove that these surfaces are counterexamples to infinitesimal Torelli. Moreover we give a complete solution for infinitesimal Torelli for elliptic surfaces (with a section)  over $\Ps^1$.

Suppose $X$ is a smooth complex algebraic variety. Then the first order deformations of $X$ are parameterized by $H^1(X,\Theta_X)$, with $\Theta_X$ the tangent bundle of $X$. The isomorphism $H^{p,q}(X)=H^q(X,\Omega^p)$ and the contraction map $\Theta_X \otimes_{\mathcal{O}_X} \Omega^p_X \ra \Omega^{p-1}_X$ give a cup product map:
\[ H^1(X,\Theta_X) \otimes H^{p,q}(X) \ra H^{p-1,q+1}(X). \]
and from this one obtains the infinitesimal period map
\[ \delta_k: H^1(X,\Theta_X) \ra \oplus_{p+q=k} \Hom (H^{p,q}(X),H^{p-1,q+1}(X)).\]

The (holomorphic) map $\delta_k$ is closely related to the period map. Assume that $\varphi: \mathcal{X} \ra \mathcal{B}$ is a proper, smooth, surjective holomorphic map between complex manifolds with connected fibers, and that for all $t\in \mathcal{B}$ the vector space $H^k(X_t,\C)$ carries a Hodge structure of weight $k$, with $X_t:=\varphi^{-1}(t)$. Fix a point $0\in \mathcal{B}$. Let $U$ be a small simply connected open neighborhood of $0$, such that the Kodaira-Spencer at the point $0$ is injective. Define
\[\mathcal{P}^{p,k} : U \ra Grass(\sum_{i\geq p} h^{i,k-i}, H^k(X_0)) \]
by sending $t$ to $(\oplus_{i\geq p} H^{i,k-i}(X_t)) \subset H^k(X_0,\C)$, via the identification $H^k(X_0,\C)\cong H^k(X_t,\C)$. (Note that $U$ is simply connected.) Then the differential of $\oplus_d \mathcal{P}^{p,k}$ is injective if and only if $\delta_k$ is injective. (See \cite[Section 2]{Sai} or \cite[Chapter 10]{Voi}.)

We say that $X$ satisfies infinitesimal Torelli if and only if $\delta_{\dim(X)}$ is injective. Note that if $X$ is an elliptic surface with a section and base $\Ps^1$ then $\delta_k$, for $k\neq \dim X$ is the zero-map.

If $X$ is a rational surface, then there is no variation of Hodge structures. If $X$ is a $K3$ surface, then infinitesimal Torelli follows from \cite{PSS}. This section will focus on the case $p_g(X)>1$.

There is an easy sufficient condition of Ki\u\i (\cite{Kii}), Lieberman, Peters and Wilsker (\cite{LPW}) for checking infinitesimal Torelli for manifolds with divisible canonical bundle. The following result is a direct consequence of {\cite[Theorem 1']{LPW}}.

\begin{Thm} \label{LPWThm} Let $X$ be a compact K\"ahler $n$-manifold, with $p_g(X)>1$. Let $\mathcal{L}$ be a line bundle such that
\begin{enumerate}
\item  $\mathcal{L}^{\otimes k}=\Omega^n_X$ for some $k>0$.
\item the linear system corresponding to $\mathcal{L}$ has no fixed components of codimension 1.
\item $H^0(X,\Omega^{n-1}_X\otimes \mathcal{L})=0. $
\end{enumerate}
Then  $\delta_n$ is injective.
\end{Thm}

We want to apply  the above theorem, when $X$ is an elliptic surface and $\mathcal{L}=\mathcal{O}_X(F)$.

\begin{Lem}\label{difflem} Let $\pi: X \rightarrow \Ps^1$ be an elliptic surface. Assume that $X$ is not birational to a product $C\times \Ps^1$. Then for $n>0$
\[ \dim H^0(X,\Omega^1(nF)) = \left\{\begin{array}{ll} n-1 & \mbox{if $j(\pi)$ is not constant} \\
 n-1 + \max(0,n+d+1)& \mbox{if $j(\pi)$ is constant,}\end{array} \right.\] 
where $d=\deg(\Li)-\# \{P\in C(\C) \mid \pi^{-1}(P) \mbox{ singular }\}$. \end{Lem}

\begin{proof} By \cite[Prop. 4.4 (I)]{Sai} we know that $\pi_* \Omega^1_X \cong \Omega^1_{\Ps^1}$ if $j(\pi)$ is not constant, which gives the first case.

If $j(\pi)$ is constant then we have the following exact sequence (by \cite[Prop. 4.4 (II)]{Sai}):
\[ 0 \ra \Omega^1_{\Ps^1} \ra \pi_* \Omega^1_X \ra \mathcal{O}_{\Ps^1}(d) \ra 0. \]
Tensoring with $\mathcal{O}_{\Ps^1}(n)$ gives
\[ \dim H^0(X,\Omega^1_X(nF))= \dim H^0(\Ps^1,\Omega^1_{\Ps^1}(n)) + \dim H^0(\Ps^1,\mathcal{O}_{\Ps^1}(n+d)) = n-1 + \max(0, n+d+1), \]
using that $\dim H^1(\Ps^1,\Omega^1_{\Ps^1}(n))=0$.
\end{proof}

\begin{Cor}\label{torcor} Suppose $\pi : X \ra \Ps^1$ is an elliptic surface (cf. Assumption~\ref{mainass}), such that $p_g(X)>1$. Suppose that $j(\pi)$ is non-constant or $\pi$ is not extremal then $X$ satisfies infinitesimal Torelli.\end{Cor}

\begin{proof} Let $F$ be a smooth fiber of $\pi$. Note that $\mathcal{O}(F)^{\otimes (p_g(X)-1)} = \Omega^2_X$ (see \cite[Theorem 2.8]{Sd}) and  $|F|$ is the elliptic fibration, hence without base-points. 

We claim that $\dim H^0(X,\Omega^1(F))=0$. Suppose if this is not  the case then by Lemma~\ref{difflem} $j(\pi)$ is constant and $\pi$ has at most $\deg(\Li)+1$ singular fibers. From Lemma~\ref{cstcrt} it follows that $j(\pi)$ is extremal.

Apply now Theorem~\ref{LPWThm} with $\mathcal{L}=\mathcal{O}(F)$.
\end{proof}

\begin{Prop}\label{KurPrp}  Let $\varphi: \mathcal{X} \ra \mathcal{B}$ be a family of surfaces with $p_g(X_0)>0$, and such that for all $t\in \mathcal{B}$ the set $\{ t' \mid X_t\cong X_{t'}\}$ is zero dimensional. Let $r$ be the Picard number of a generic member of the family $\varphi$.
Suppose that if $p_g(X_0)>1$ and
\[ \dim \mathcal{B} \geq \frac{1}{2} p_g(X) (h^{1,1}-r). \]
or $p_g(X_0)=1$ and
\[ \dim \mathcal{B} \geq  (h^{1,1}-r). \]
Then for all $t$, the surface $X_t$ does not satisfy infinitesimal Torelli. \end{Prop}

\begin{Cor} Let $\psi: \mathcal{X}\ra \mathcal{B}$ be a non-tivial family of surfaces such that $\rho(X_t)=h^{1,1}(X_t)$ for all $t$. Then $X_t$ does not satisfy infinitesimal Torelli.\end{Cor}

\begin{Cor}\label{antitorcor} Let $\pi: X \ra \Ps^1$ be an extremal elliptic surface with constant $j$-invariant and $p_g(X)>1$. Then $X$ does not satisfy infinitesimal Torelli.
\end{Cor}

\begin{proof}From Remark~\ref{modrem} it follows that $X$ is a member of positive dimensional family of surfaces with $r=h^{1,1}(X_0)$.
\end{proof}

\begin{proof}[Proof of Proposition~\ref{KurPrp}]

 Fix a base point $0 \in \mathcal{B}$ and let $X$ be isomorphic to the fiber over $0$. For any $t\in \mathcal{B}$ denote by $X_t$ the fiber over $t$. Let $\Lambda$ be a lattice of rank $r$, such that $\Lambda \hookrightarrow NS(X_0)$, and we can fix the identification of Hodge structures $H^2(X_0,\C)\cong H^2(X_t,\C)$ such that $\Lambda \hookrightarrow NS(X_t)$. 

Let $T(X_t)$ be the orthogonal complement of the image of $\Lambda$ in $H^2(X_t,\Z)$. Then $T(X_t)\otimes \C$ carries a sub-Hodge structure, and we consider  variation of this Hodge structure. (cf. \cite[Section 6]{vGT}). Since the variation of the Hodge structure on $H^2(X,\C)$ is determined by the variation of the sub-Hodge structures $\Lambda$ and $T$, and $\Lambda$ remains of pure type $(1,1)$
we have the following diagram

\[ \begin{array}{ccccc}
\mathcal{X} & \longrightarrow &\mathcal{X}_{k} \\
\downarrow & & \downarrow \\
\mathcal{B} & \stackrel{\psi}{\longrightarrow} & \mathcal{B}_k \\
  & & \downarrow \rho \\
    &&   H^1(X_0,\Theta_{X_0}) & \stackrel{\delta_2}{\longrightarrow} & \Hom(H^{2,0}(X_0),H^{1,1}(X_0)) \oplus \Hom (H^{1,1}(X_0),H^{0,2}(X_0)) \\

&& \uparrow && \uparrow\\
&& Im(\rho \psi) & \longrightarrow &\Hom(T^{2,0}(X_0),T^{1,1}(X_0)) \oplus \Hom (T^{1,1}(X_0),T^{0,2}(X_0)),\\
\end{array} \]
with $\mathcal{X}_k$  the Kuranishi family of $X_0$ and $\mathcal{B}_k$ its base.

From this diagram one deduces $\dim Im(\rho\psi)\leq 2(h^{1,1}-r)h^{2,0}$. Using Serre duality and a result of Griffiths one can show (cf. \cite[\S 6]{vGT}) that $Im(\rho \psi)$ has dimension at most $\frac{1}{2} (h^{1,1}-r)h^{2,0}$, if $h^{2,0}>1$ and $Im(\rho \psi)$ has dimension at most $h^{1,1}-r$ if $h^{2,0}=1$.
\end{proof}

\begin{Rem} The referee pointed out that the  factorization
\[ \mathcal{B} \ra \Hom(T^{2,0}(X_0),T^{1,1}(X_0)) \ra  \Hom(H^{2,0}(X_0),H^{1,1}(X_0)) \]
can be proven as follows. 

Let $t\in B_k$ be  a direction contained for which the image of $\Lambda$ remain divisors. Then the kernel of cup product with $t$ consists of the classes of type $(1,1)$ infinitesimally fixed in that direction.
\end{Rem}

\begin{Thm}\label{TorThm} Let $\pi: X \ra \Ps^1$ be an elliptic surface (cf. Assumption~\ref{mainass}). Assume that $p_g(X)>1$. Then $X$ does {\em  not} satisfy infinitesimal Torelli if and only if $j(\pi)$ is constant and $\pi$ is extremal. \end{Thm}

\begin{proof} Combine Corollary~\ref{torcor} and Corollary~\ref{antitorcor}.
\end{proof}



\begin{Rem} Note that the hyperelliptic surfaces form a family of elliptic surfaces with $p_g(X)=0$, so they do not satisfy infinitesimal Torelli.\end{Rem}

\begin{Rem}  Chakiris (\cite[\S 4]{Cha}) gave different formulae for $\dim H^0(X,\Omega^1(nF))$. He uses them to deduce a formula for $\dim H^1(X,\Theta_X)$, which he uses to prove that generic global Torelli holds. Even with the use of these incorrect formulae his proof of generic global Torelli seems to remain valid, after a small modification. His formulae would imply that infinitesimal Torelli holds for any elliptic surface with a section. The same erroneous formulae leads Beauville (\cite[p. 13]{Beau}) to state in a survey paper on Torelli problems that infinitesimal Torelli holds for {\em arbitrary}  elliptic surface with a section. Theorem~\ref{TorThm} shows instead that this is true under the condition that $j(\pi)$ is not constant or $\pi$ is not extremal.
\end{Rem}

The argument used to prove that several elliptic surfaces satisfy infinitesimal Torelli, relies heavily on $C=\Ps^1$. Sait\=o proved using other techniques that if $\pi: X \ra C$ is an elliptic surface with non-constant $j$-invariant then $X$ satisfy infinitesimal Torelli. We consider now the case that the $j$-invariant is constant and we try to find surfaces for which infinitesimal Torelli does not hold.

\begin{Lem}  Let $\varphi:\mathcal{X}\ra\mathcal{B}$ be a  family of elliptic surfaces with $p_g(X_0)>1$, constant $j$-invariant and $s$ singular fibers. Let $g(C)$ be the genus of the base curve of a generic member of this family. Suppose that for all $t$ we have that $\{t'\mid X_t \cong X_{t'}\}$ is zero-dimensional and 
\[ \dim \mathcal{B} > (s-\deg(\Li)+g(C)-1)h^{2,0}.\]
Then for all $t$, $X_t$ does not satisfy infinitesimal Torelli.
\end{Lem}

\begin{proof}

Note that $\rho(X_t) \geq \rho_{tr}(\pi)$ for all $t\in \mathcal{B}$ and 
\[h^{1,1}(X)-\rho_{tr}(\pi) = 2s-2\deg(\Li) - 2+2g(C). \]
Apply now Proposition~\ref{KurPrp}.
\end{proof}

\begin{Exa} Let $\psi:\mathcal{X} \ra \mathcal{B}$ be a maximal-dimensional family of elliptic surfaces with $2\deg(\Li)$ singular fibers of type $I_0^*$. Then $\dim \mathcal{B}=3g-3+2\deg(\Li)+1$. 

We have that $h^{2,0}>1$ and 
\[ \dim\mathcal{B}> (s-\deg(\Li)+g(C)-1)h^{2,0} \]
holds if and only if $h^{2,0}=2, g(C)\in\{ 1,2,3\}$ or $h^{2,0}=3, g(C)=4$. In all these cases any member of the family $\psi$ is a counterexample to Infinitesimal Torelli. \end{Exa}

\begin{Exa} Let $\psi:\mathcal{X} \ra \mathcal{B}$ be a maximal family of elliptic surfaces with $j$-invariant 0 or 1728 and $s$ singular fibers. Then $\dim \mathcal{B} = 3g-3+s$. 

Using $h^{2,0}=\deg(\Li)+g(C)-1$, we obtain that 
\[ \dim\mathcal{B}> (s-\deg(\Li)+g(C)-1)h^{2,0}\]
holds if and only if
\[ s < \frac{\deg(\Li)^2-g(C)^2+5g(C)-4)}{\deg(\Li)+g(C)-2}.\]
From Noether's condition, the smallest $s$ that is possible is $\frac{6}{5}\deg \Li$. From this one deduces
\[ \deg(\Li) < -3g(C)+6+\sqrt{4g(C)^2-11g(C)+16}.\]

All combinations of $g(X)$ and $p_g(X)$ satisfying the above conditions are mentioned in the table below. There $s_{max}$ denotes the maximum number of singular fibers such that the maximal family with $s_{max}$ singular fibers has $\dim(\mathcal{B})$ is larger then the upper bound for the dimension of the period domain. The columns with $[ 6/5\deg(\Li) ]$ and $[4/3\deg(\Li)]$ denote the minimal number of singular fibers for an elliptic surface with $j$-invariant 0 or 1728, whenever this is smaller then $s_{max}$.
\[\begin{array}{|c|c|c|c|c|c|}
\hline
g(C) & p_g(X)  &\deg(\Li) &s_{max}& [ 6/5\deg(\Li) ] & [4/3\deg(\Li)]\\
\hline
1    &   2 &    2   &3& 3  & 3\\
1    & 3 & 3 & 4&4 & 4\\
1   &  4 & 4  & 5&5 & -\\
1   & 5 & 5 & 6&6&  - \\
2   & 2  & 1 & 2&2& 2 \\ 
\hline
\end{array}\]
\end{Exa}

\section{Twisting}\label{twist}
In this section we study the behavior of $h^{1,1}(X)-\rho_{tr}(X)$ under twisting, when $X$ is a Jacobian elliptic surface.

To a Jacobian elliptic surface  $\pi: X \ra C$, we can associate an elliptic curve in $\Ps^2_{\C(C)}$  corresponding to the generic fiber of $\pi$. Denote this elliptic curve by $E_1$. Vice versa to an elliptic curve $E'/\C(C)$ we can associate an elliptic surface $\pi': X' \ra C$.

Let $\pi_i : X_i \ra C, i=1,2$ be two Jacobian elliptic surfaces, such that $j(\pi_i)\neq 0,1728$. Then there exists an isomorphism $\varphi: X_1\ra X_2$ such that $\pi_2\varphi=\pi_1$ if and only if the associated elliptic curves $E_1$ and $E_2$ are isomorphic over $\C(C)$. 

The last statement is equivalent to $j(E_1)=j(E_2)$ and the quotients of the minimal discriminants of $E_1/\C(C)$ and $E_2/\C(C)$ is a 12-th power (in $\C(C)^*$).

If $j(E_1)=j(E_2)$ then under our assumptions the quotient of the minimal discriminants is a 6-th power, say $u^6$. From this it follows that $E_1$ and $E_2$ are isomorphic over $\C(C)(\sqrt{u})$. We call $E_2$ the twist of $E_1$ by $u$ and denote this by $E_1^{(u)}$. Actually, we are not interested in the function $u$, but in the places at which the valuation of $u$ is odd.

\begin{Def} Let $\pi:X \ra C$ be a Jacobian elliptic surface. Fix $2n$ points $P_i \in C(\C)$. Let $E/\C(C)$ be the Weierstrass model of the generic fiber of $\pi$.

A Jacobian elliptic surface $\pi': X' \ra C$ is called a \emph{quadratic twist} of $\pi$ by $(P_1,\ldots,P_n)$ if the Weierstrass model of the generic fiber of $\pi'$ is isomorphic to $E^{(f)}$, where $E^{(f)}$ denotes the quadratic twist of $E$ by $f$ in the above mentioned sense 
and $f\in\C(C)$ is a function such that  $v_{P_i}(f) \equiv 1 \bmod 2$ and $v_Q(f) \equiv 0 \bmod 2$ for all $Q\not \in \{P_i\}$.

  A \emph{*-minimal twist} of $\pi$ is a twist $\tilde{\pi}:\tilde{X} \rightarrow C$ such that none of the fibers are of type $II^*, III^*, IV^*$ or $I_\nu^*$ and at most 1 fiber is of type $I_0^*$.
\end{Def}

The existence of the twist follows easily from the fact that $Pic^0(C)$ is divisible. If $g(C)>0$ and we fix $2n$ points $P_1,\ldots P_{2n}$ then there exist $2^{2g(C)}$ twists by $(P_1,\ldots ,P_{2n})$.

Note that a *-minimal twist is a twist for which $p_g(X)$ (and $p_a(X)$) is minimal. Later on we will introduce another notion of minimality: the twist for which $h^{1,1}(X)-\rho_{tr}(\pi)$ is minimal.

A *-minimal twist of an elliptic curve need not be unique, but the configuration of the singular fibers of any two *-minimal twists of the same surface are equal.

If $P$ is one of the $2n$ distinguished points, then the fiber of $P$ changes in the following way (see \cite[V.4]{Mi1}).
\[
I_\nu \leftrightarrow I^*_\nu (\nu \geq 0) \;\;\; \;\;
II \leftrightarrow    IV^* \;\;\;\;\;
III  \leftrightarrow  III^* \;\;\;\;\;
IV  \leftrightarrow  II^*
\]

\begin{Lem}\label{twstcon} Let $\pi: X \ra C$ be a Jacobian elliptic surface. Let $P_i\in C(\C)$ be $2n$ points. Let $\pi':X' \ra C$ be a twist by $(P_i)$. Then
\[ h^{1,1}(X')-\rho_{tr}(\pi')= h^{1,1}(X)-\rho_{tr}(\pi) + \sum_{i=1}^{2n} c_{P_i} \]
with 
\[c_{P_i}=\left\{ \begin{array}{cl} 
1 & \mbox{if $\pi^{-1}(P_i)$ is of type $I_0,IV^*, III^*$ or $II^*$,}\\
0 & \mbox{if $\pi^{-1}(P_i)$ is of type $I_\nu$ or $I_\nu^*$, with $\nu>0$,}\\
-1 & \mbox{if $\pi^{-1}(P_i)$ is of type $II, III, IV,$ or $I_0^*$.}\\ \end{array} \right. \]
\end{Lem}

\begin{proof} Suppose $\pi^{-1}(P_i)$ is of type $I_0$. Then $\pi'^{-1}(P_i)$ is of type $I_0^*$. The Euler characteristic of this fiber is 6, so this point causes  $h^{1,1}$ to increase by 5. An $I_0^*$ fiber has 4 components not intersecting the zero-section. Hence $\rho_{tr}$ increases by 4.

The other fiber types can be done similarly.
\end{proof}

\begin{Lem}\label{mincar} Given a Jacobian elliptic surface $\pi:X \rightarrow C$ with non-constant $j$-invariant. There exist finitely many twists $\pi'$ of $\pi$ such that the non-negative integer $h^{1,1}-\rho_{tr}$ is minimal under twisting. These twists are characterized by $b(\pi')=c(\pi')=0$. \end{Lem}

\begin{proof} It is easy to see that there are at most finitely many twists with $c(\pi')=0$. Hence it suffices to show that $b(\pi')=c(\pi')=0$ if $h^{1,1}(X')-\rho_{tr}(\pi')$ is minimal under twisting. From Lemma~\ref{twstcon} it follows that it suffices to show that for any elliptic surface there exists a twist with $b(\pi')=c(\pi')=0$.

  Consider a *-minimal twist $\tilde{\pi}:\tilde{X} \rightarrow C$. Note that $e(\tilde{\pi})>0$ (otherwise the $j$-invariant would be constant.)

Suppose $b(\tilde{\pi})+c(\tilde{\pi})$ is even. Twist by all points with a fiber of type $II, III, IV, I_0^*$. The new elliptic surface has $b=c=d=0$.

Suppose $b(\tilde{\pi})+c(\tilde{\pi})$ is odd. Twist by all points with a fibers of type $II, III, IV, I_0^*$ and one point with a fiber of type $I_\nu$. The new elliptic surface has $b=c=0, d=1$.
\end{proof}

\begin{Rem}\label{unique} The classification (in \cite{SZ}) of extremal $K3$ surfaces is a classification of the root lattices corresponding to the singular fibers. In general one {\em cannot} decide which singular fibers correspond to these lattices, since each of the pairs $(I_1,II)$, $(I_2,III)$ and $(I_3,IV)$ give rise to the same lattice ($A_0$,$A_1$, $A_2$). From the lemma above it follows that this is not a problem when $\pi$ is extremal.\end{Rem}

\begin{Prop} Let $\pi_1: X_1 \rightarrow C$ be an elliptic surface with $j(\pi_1)$ non-constant. Let $\pi: X\rightarrow C$ be a *-minimal twist of $Jac(\pi_1)$ with associated line bundle $\Li$. Let $\tilde{\pi}: \tilde{X} \rightarrow C$ be a twist for which $h^{1,1}-\rho_{tr}$ is minimal. Then
  \[ h^{1,1}(\tilde{X})-\rho_{tr}(\tilde{X})=2g(C)-2\deg(\Li)-2+\# \{\mbox{singular fibers for } \pi \}. \]
\end{Prop}

\begin{proof} From Proposition~\ref{fibcrit} and the Lemmas~\ref{mincar} and~\ref{twstcon} we have $\deg(\tilde{\Li})=\deg(\Li)+(d(\tilde{\pi})+a(\tilde{\pi})-c(\pi))/2$, $d(\tilde{\pi})+e(\tilde{\pi})=e(\pi)$, $a(\tilde{\pi})=b(\pi)$. This yields
 \begin{eqnarray*}
   h^{1,1}(\tilde{X})-\rho_{tr}(\tilde{X})& = &2g(C)-2 \deg(\tilde{\Li})-2+2(a(\tilde{\pi})+d(\tilde{\pi}))+e(\tilde{\pi}) \\
    & = & 2g(C)-2\deg(\Li) -2 + b(\pi) +c(\pi) +e(\pi).
  \end{eqnarray*}
Finally note that $a(\pi)=d(\pi)=0$.
\end{proof}

\begin{Cor} \label{ExtCar} Let $\tilde{\pi}: \tilde{X} \rightarrow C$ be an elliptic surface with $j(\pi)$ non-constant, then $\tilde{\pi}$ is extremal if and only if $\tilde{\pi}$ has no fibers of type $II, III, IV$ or $I_0^*$ and the *-minimal twist of its Jacobian $\pi:X \rightarrow C$ has $2\deg(\Li)+2-2g(C)$ singular fibers.
\end{Cor}

\section{Configurations of singular fibers}\label{singfib}
In order to apply the results of the previous section, we need to know which elliptic surfaces have a minimal twist with $2\deg(\Li)+2-2g(C)$ singular fibers.

We need the following definition: 

\begin{Def} A function $f: C \rightarrow \Ps^1$ is called of $(3,2)$-type if the ramification indices of the points in the fiber of $0$ are at most $3$, and in the fiber of $1728$ are at most $2$.\end{Def}

For the connection between functions of $(3,2)$-type and certain subgroups of $SL_2(\Z)$ see \cite{Beaufour}.

\begin{Prop}\label{NumFib} Let $\pi: X \rightarrow C$ be an elliptic surface with $j(\pi)$ non-constant, such that $\pi$ is a *-minimal twist. Then $j(\pi)$ is of $(3,2)$-type and not ramified outside $0, 1728, \infty$ if and only if there are $2\deg(\Li)+2-2g(C)$ singular fibers. \end{Prop}

\begin{proof}
  Denote by
  \begin{itemize}
  \item $n_2$ the number of fibers of $\pi$ of type $II$.
  \item $n_3$ the number of fibers of $\pi$ of type $III$.
  \item $n_4$ the number of fibers of $\pi$ of type $IV$.
  \item $n_6$ the number of fibers of $\pi$ of type $I_0^*$.
  \item $m_\nu$ the number of fibers of $\pi$ of type $I_\nu$. ($\nu>0$.)
  \end{itemize}
Let $r=\sum \nu m_\nu$.

The ramification of $j(\pi)$ is as follows (using \cite[Lemma IV.4.1]{Mi1}).

Above $0$ we have $n_2$ points with  ramification index 1 modulo 3, we have $n_4$ point with index 2 modulo 3 and at most $(r-n_2-2n_4)/3$ points with index 0 modulo 3. In total, we have at most $n_2+n_4+(r-n_2-2n_4)/3$ points in $j(\pi)^{-1}(0)$.

Above $1728$ we have $n_3$ points with index 1 modulo 2 and at most $(r-n_3)/2$ points with  index 0 modulo 2. So $j(\pi)^{-1}(1728)$ has at most $n_3+(r-n_3)/2$ points.

Above $\infty$ we have $\sum m_\nu$ points.

Collecting the above gives
\[ \#j(\pi)^{-1}(0) + \#j(\pi)^{-1}(1728) + \#j(\pi)^{-1}(\infty) \leq 2/3 n_2+ 1/2 n_3+ 1/3 n_4 + 5/6r + \sum m_\nu \]
with equality if and only if $j(\pi)$ is of $(3,2)$-type.

Hurwitz' formula implies that
\[ r+2-2g(C) \leq \#j(\pi)^{-1}(0) + \#j(\pi)^{-1}(1728) + \#j(\pi)^{-1}(\infty) \]
with equality if and only if there is no ramification outside $0, 1728$ and $\infty$.

So
\[r\leq 12g(C)-12+4 n_2+3 n_3+2 n_4+6\sum m_\nu \]
and by Proposition~\ref{Noether}
\[ \sum i n_i + r =12\deg(\Li). \]

Substituting gives
\[ 2\deg(\Li) +2-2g(C)\leq n_2+n_3+n_4+n_6+\sum b_\nu\]
with equality if and only if $j(\pi)$ is not ramified outside $0, 1728$ and $\infty$ and $j(\pi)$ is of $(3,2)$-type. \end{proof}

This enables us to prove
\begin{Thm} \label{MainThm} Suppose $\pi:X\rightarrow C$ is an elliptic surface with non-constant $j$-invariant, then the following three are equivalent
\begin{enumerate}
\item $\pi$ is extremal 
\item $j(\pi)$ is of $(3,2)$-type, not ramified outside $0, 1728, \infty$  and $\pi$ has no fibers of type $II, III, IV$ or $I_0^*$.
\item The minimal twist $\pi'$ of $Jac(\pi)$ has $2\deg(\Li)+2-2g(C)$. 
\end{enumerate}
\end{Thm}

\begin{proof} Apply Proposition~\ref{NumFib} to Corollary~\ref{ExtCar}.\end{proof}

\begin{Rem} Frederic Mangolte pointed out to me that the equivalence of (1) and (2) was already proved in \cite{Nori}.\end{Rem}

\begin{Rem} Consider functions $f:C \rightarrow \Ps^1$ up to automorphisms of $C$. If we fix the ramification indices above $0, 1728, \infty$ and demand that $f$ is unramified at any other point, then there are only finitely many $f$ with that property. A small deformation of $f$ in $Mor_d(C, \Ps^1)$, the moduli space of morphisms $C \rightarrow \Ps^1$ of degree $d$, has more critical values.

  So the $j$-invariants of extremal elliptic surface lie ``discretely'' in $Mor_d(C,\Ps^1)$. By Lemma~\ref{mincar}, to any $j$-invariant there correspond only finitely many extremal elliptic surfaces (by Lemma~\ref{mincar}). So extremal elliptic surfaces over $\Ps^1$ with geometric genus $n$ and non-constant $j$-invariant, lie discretely in the moduli space of elliptic surfaces over $\Ps^1$ with geometric genus $n$.
  \end{Rem}

Suppose that $\pi: X \rightarrow \Ps^1$ has $2p_g(X)+4$ singular fibers of type
$I_\nu$ and no other singular fibers. Let $f: \Ps^1 \rightarrow \Ps^1$ be a
cyclic morphism ramified at two points $P$ such that $\pi^{-1}(P)$ is singular.

Fastenberg (\cite[Theorem 2.1]{Fa}) proved that then the base changed surface has
Mordell-Weil rank 0. In fact, she proved that the base-changed surface is extremal. The
first surface is also extremal (by Proposition~\ref{ExtCar}). A slightly more
general variant is the following.

\begin{Exa}\label{FastExa} Suppose $\pi: X \rightarrow C$ is an extremal elliptic surface. Let $f:C' \rightarrow C$ be a finite morphism.

 Then the base-changed elliptic surface $\pi': X' \rightarrow C'$ is extremal if $f$ is not ramified outside the set of points $P$, such that $\pi^{-1}(P)$ is
multiplicative or potential multiplicative.  

In that case the composition $j':
C'\rightarrow C \stackrel{j}{\rightarrow}\Ps^1$ is not ramified outside $0, 1728$
and $\infty$ and the ramification indices above $0$ and $1728$ are at most 3 and
2. Moreover there are no fibers of type $II, III, IV$ or $I_0^*$.
\end{Exa}

An easy calculation shows that all elliptic surfaces for which Fastenberg's results (\cite[Theorem 1]{Fa}) hold, are either extremal elliptic surfaces or have a twist which is extremal. Moreover these surfaces have no fibers of type $I_0^*$, hence all elliptic surfaces for which her results hold lie discretely in the moduli spaces mentioned above.
  \section{Mordell-Weil groups of extremal elliptic surfaces}\label{MW}
It remains to classify which Mordell-Weil groups can occur. To this we can give only a partial answer to this.

Let $m,n \in \Z_{\geq 1}$ be such that $m|n$ and $n>1$. Recall from section~\ref{Intro} that $X_m(n)$ is the modular curve parameterizing triples $((E,O),P,Q)$, such that $(E,O)$ is an elliptic curve, $P\in E$ is a point of order $m$ and $Q\in E$ a point of order $n$.

If $(m,n)\not \in \{ (1,2), (2,2), (1,3), (1,4), (2,4) \}$ then there exists a universal family for $X_m(n)$, which we denote by $E_m(n)$. Denote by $j_{m,n}: X_m(n) \rightarrow \Ps^1$ the map usually called $j$.

From the results of \cite[\S 4 \& \S 5]{ShEMS} it follows that $E_m(n)$ is an extremal elliptic surface. The following theorem explains how to construct many examples of extremal elliptic surfaces with a given torsion group.

\begin{Thm} Fix $m,n\in \Z_{\geq 1}$ such that $m|n$ and $(m,n)\not \in \{(1,1), (1,2), (2,2), (1,3), (1,4),(2,4) \}$. Let $j\in \C(C)$ be non-constant then there exists a unique elliptic surface $\pi:X \rightarrow C$, with $j(\pi)=j$ and $MW(\pi)$ has $\Z/n\Z\times \Z/m\Z$ as a subgroup if and only if $j$ is of $(3,2)$-type, not ramified outside $0,1728, \infty$ and $j=j_{n,m} \circ g$ for some $g: C \rightarrow X_m(n)$.

The extremal elliptic surface with $j(\pi)=j$ and $\Z/n\Z \times \Z/m\Z$ is a subgroup of $MW(\pi)$ is the unique surface with only singular fibers of type $I_\nu$.
\end{Thm}

\begin{proof}  Let $\pi: X \rightarrow C$ be an elliptic surface, such that $MW(\pi)$ has a subgroup isomorphic to $G:=\Z/m\Z \times \Z/n\Z$. Then $j:C\rightarrow \Ps^1$ can be decomposed in $g: C \rightarrow X_m(n)$ and $j_{m,n}: X_m(n) \rightarrow \Ps^1$, and $X$ is isomorphic to $E_m(n) \times_{X_m(n)} C$.

  Conversely, for any base change $\pi'$ of $\varphi_{m,n}: E_m(n) \rightarrow X_m(n)$, the group $MW(\pi')$ has $G$ as a subgroup.

  Moreover since $\varphi_{m,n}$ has only singular fibers of type $I_\nu$, the same holds for $\pi'$. An application of Theorem~\ref{MainThm} concludes the proof.
\end{proof}

\begin{Rem} Let $n=2$ and $m\leq 2$. Then any elliptic surface with such that $j(\pi)=j_{m,n}\circ g$, has $G=\Z/n\Z \times \Z/m\Z$ as a subgroup of $MW(\pi)$.\end{Rem}

\section{Uniqueness}\label{uni}
Artal Bartolo, Tokunaga and Zhang (\cite{ATZ}) expect that an extremal elliptic surface is determined by
the configuration of the singular fibers and the Mordell-Weil group. More precise, they raise
the following question:
\begin{Que} Suppose $\pi_1: X_1 \rightarrow \Ps^1$ and $\pi_2: X_2 \rightarrow \Ps^1$ are extremal semi-stable elliptic surfaces, such that $MW(\pi_1)\cong MW(\pi_2)$ and the configurations of singular fibers of $\pi_1$ and $\pi_2$ are the same.

  Are $X_1$ and $X_2$ isomorphic, and if so, is there then an isomorphism that respects the fibration and the zero section?
\end{Que}

By \cite[Theorem 5.4]{MP} this is true in the case where $X_1$ and $X_2$ are rational elliptic surfaces.

In the case where $X_1$ and $X_2$ are $K3$ surfaces the answer is the following theorem.

\begin{Thm} There exists precisely 19 pairs $(\pi_1: X_1 \rightarrow \Ps^1, \pi_2: X_2 \rightarrow \Ps^1)$ of extremal elliptic $K3$ surface, such that $\pi_1$ and $\pi_2$ have the same configuration of singular fibers, $MW(\pi_1)$ and $MW(\pi_2)$ are trivial, and $X_1$ and $X_2$ are not isomorphic. Of these pairs 13 are semi-stable. There is an unique pair for which the above holds with $MW(\pi_1)=MW(\pi_2)=\Z/2\Z$, which is not semi-stable.\end{Thm}

\begin{proof}  From \cite[Table 2]{SZ} there exist 19 pairs of surfaces $(X_1,X_2)$ such that the transcendental lattices of $X_1$ and $X_2$ lie  in  distinct $SL_2(\Z)$-orbits,  they admit elliptic fibrations $\pi_i: X_i \rightarrow \Ps^1$  such that $MW(\pi_1)=MW(\pi_2)=0$ and the contributions of the singular fibers as sub-lattices of the N\'eron-Severi lattice coincide.

  From Remark~\ref{unique} we know that for extremal elliptic surfaces the sub-lattices determine the singular fibers. Since the transcendental lattices are in different $SL_2(\Z)$-orbits, the surfaces are not isomorphic.

  The rest of the statement follows from the same Table.
\end{proof}

\begin{Rem}
From these surfaces one should be able to construct other pairs of extremal elliptic surfaces with isomorphic Mordell Weil groups, and the same configuration of singular fibers, such that the geometric genus is higher then 1.

Start with two non-isomorphic extremal elliptic $K3$ surfaces
with the same configuration of singular fibers and the same Mordell Weil group.
Then the $j$-invariant of both surfaces are unequal modulo automorphism of
$\Ps^1$.

We can base-change both surfaces in such a way that the base-changed surfaces
remain extremal (cf. Example \ref{FastExa}), they have the same configuration of singular fibers and their
$j$-invariants are unequal modulo an automorphism of $\Ps^1$.

The configuration of singular fibers gives restrictions on the possibilities for the torsion part of the Mordell-Weil group. One can hope that this is sufficient to prove that the Mordell-Weil groups are isomorphic.

Note that the base-changed surfaces are not isomorphic, since a surface which is not
a $K3$ surface has at most one elliptic fibration.
\end{Rem}

\begin{Rem} (\textbf{``case 49''}) In \cite{ATZ} it is proven that there are two elliptic surfaces with $MW(\pi_i)=\Z/5\Z$ and singular fibers $2I_1, I_2, 2I_5, I_{10}$ and there exists no isomorphism between them that respects the fibration. In \cite{SZ} it is proven that both surfaces are isomorphic. For any other pair of (semi-stable) extremal elliptic $K3$ surfaces with $\# MW(\pi)>4$ and the same singular fibers configuration, there exists an isomorphism which respects the fibration. (\cite[Theorem 0.4]{ATZ})
\end{Rem}

\section{Extremal elliptic surfaces with $p_g=1, q=1$}\label{cla}
An elliptic surface with a section and $q=1$ needs to have a genus 1 base curve.
This implies that for an extremal elliptic surface with $p_g=1, q=1$ we have
$\deg(\Li)=1$.

The minimal twist of an extremal elliptic surface with $p_g=1,q=1$ has two singular fibers.

All possible pairs of fiber types such that the sum of the Euler characteristics is 12, are given in the following table.
\[ \begin{array}{|c|c|c|c|c|c|c|c|c|c|c|}\hline
F_1&I_{11} & I_{10}& I_9& I_8& I_7 & I_6 &  I_{10}& I_9& I_8& I_6\\
F_2& I_1& I_2& I_3& I_4& I_5& I_6& II& III& IV& I_0^*\\ 
\hline
   & & \mbox{\cite{I10I2}}&\mbox{\cite{I9I3}}&\mbox{\cite{I8I4}}&&\mbox{\cite{I6I6}}&&&&\\
\hline
\end{array}\]
Several of these surfaces are already described in the literature. See the above mentioned references.

\begin{Rem} Note that the 6 configurations with two $I_\nu$ fibers are already extremal. The configurations with one additive and one multiplicative fibers are not extremal. If we twist by the two points with a singular fiber we obtain an extremal elliptic surface, but then the degree of $\Li$ is 2, except for the $I_6 \; I_0^*$, its extremal twist has one singular fiber of type $I_6^*$.\end{Rem}

\begin{Prop} All these configurations exist, except $I_7 \; I_5$.\end{Prop}

\begin{proof} This follows from the following lemmas. Note that the existence of elliptic surfaces over $\Ps^1$ with the below mentioned singular fibers follows from \cite{Per}.

\begin{Lem} The configurations with $I_k \; I_{12-k}$ exist for $k=2,4,6$.\end{Lem}

\begin{proof} Let $\pi:X\rightarrow \Ps^1$ be an elliptic surface with two $III$ fibers, a fiber of type $I_{k/2}$ and a fiber of type $I_{6-k/2}$, and no other singular fibers. Let $\varphi: C \rightarrow \Ps^1$ be a degree two cover ramified at the four points where the fiber of $\pi$ is singular. Then $\pi': X \times_{\Ps^1} C\ra C$ has two fibers of type $I_0^*$, a fiber of type $I_k$ and a fiber of type $I_{12-k}$. Twisting by the two points with $I_0^*$ fibers gives the desired configuration.\end{proof}

\begin{Lem} The configurations $I_8 \;IV$ and $I_6 \;I_0^*$ exist.\end{Lem}

\begin{proof} For the first let $\pi:X\rightarrow \Ps^1$ be an elliptic surface with two $III$ fibers, a fiber of type $II$ and a fiber of type $I_4$, and no other singular fibers. Let $\varphi: C \rightarrow \Ps^1$ be a degree two cover ramified at the four points where the fiber of $\pi$ is singular. Then $\pi': X \times_{\Ps^1} C\ra C$ has two fibers of type $I_0^*$ a fiber of type $I_8$ and a fiber of type $IV$. Twisting by the two points with a $I_0^*$ fibers gives the desired configuration.

For the second let $\pi:X\rightarrow \Ps^1$ be an elliptic surface with three $III$ fibers and a fiber of type $I_3$, and no other singular fibers. Let $\varphi: C \rightarrow \Ps^1$ be a degree two cover ramified at the four points where the fiber of $\pi$ is singular. Then $\pi': X \times_{\Ps^1} C\ra C$ has three fibers of type $I_0^*$ and a fiber of type $I_6$. Twisting by two points with a $I_0^*$ fiber gives the desired configuration.
\end{proof}

For the other four configurations we have a different strategy. We simply show that the $j$-map with the right ramification indices exists. This is equivalent to give the monodromy representation.

\begin{Lem} The configurations $I_{11} \;I_1$, $I_9\; I_3$, $II\; I_{10}$, $III\; I_9$ exist.\end{Lem}

\begin{proof} For the two configurations of type $I_\mu \; I_\nu$  we need to find curves $C$ and  functions $j:C \ra \Ps^1$ of degree 12 such that above $\infty$ there are two point with ramification indices $\mu$ and $\nu$, all points above $0$ have ramification index 3 and all points above $1728$ have ramification index 2. (see \cite[Lemma IV.4.1]{Mi1}.)

By the Riemann existence theorem it suffices to give two permutations $\sigma_0,\sigma_1$ in $S_{12}$, such that $\sigma_0$ is the product of 6 disjoint 2-cycles, $\sigma_1$ the product of 4 disjoint 3-cycles, and $\sigma_0\sigma_1$ is a product of a $\mu$-cycle and a $\nu$ cycle, and the subgroup generated by $\sigma_0$ and $\sigma_1$ is transitive. (See \cite[Corollary 4.10]{MiRS}.)

For $I_1 \; I_{11}$, we use
\[ (1\; 2 \; 3)(4\; 5 \; 6)(7 \; 8 \;9)(10\; 11\; 12)*(1\; 3)(2\; 4)(5\; 7)(8\; 10)(9\; 11)(6\; 12)=(1)(3\; 2 \; 5 \; 8 \; 11 \;7 \; 6 \;10 \; 9 \; 12 \; 4)\]
For $I_3 \; I_9$, we use
\[(1\; 2 \; 3)(4\; 5 \; 6)(7 \; 8 \;9)(10\; 11\; 12)*(1 \;6)(4\;9)(3\; 7)(2\; 10)(5\; 12)(8\;11)= (1 \; 4 \;7)(2\; 11\;9\;5 \; 10 \; 3 \; 8 \; 12\; 6)\]
Similarly the existence of $II \; I_{10}$, follows from
\[ (1\; 2)(3\; 4)(5\; 6)(7\; 8)(9\; 10) * (1 \; 4 \; 7)(2 \; 5 \; 8)(3 \; 6 \;9)=(1 \; 3\; 5\;7 \; 2 \; 6 \; 10 \;9 \;4 \; 8)\]
and the existence of $ III\; I_9$, follows from
\[ (2\; 3)(4 \; 5)(6\; 7)(8 \;9) * (1\; 4 \;7)(2\; 5 \; 8)(3\; 6 \;9)=(1\;5\;9\;2\;4\;6\;8\;3\;7). \]
\end{proof}

\begin{Lem} The configuration $I_7 \; I_5$ does not exists.\end{Lem}

\begin{proof}A computer search learned us that the permutations needed for the existence of $I_7\; I_5$ do not exist.\end{proof}\end{proof}

\begin{Cor} Let $k_i$ be positive integers such that $\sum k_i=12$, with $i\geq 2$, and if $i=2$ then $(k_1,k_2)\neq (7,5)$ or $(5,7)$. Then there exist a curve $C$ of genus 1, and an elliptic surface $\pi: X \rightarrow C$ such that the configuration of singular fibers of $\pi$ is $\sum I_{k_i}$.
\end{Cor}

\begin{proof} Use the monodromy representation as in \cite[Remark after Corollary 3.5]{Mi2}.
\end{proof}

\section{Elliptic surfaces with one singular fiber}\label{one}
In the previous section we proved that there exists an elliptic fibration with an $I_6$ and an $I_0^*$ fiber. Twisting by the points with a singular fiber yields an elliptic surface with one singular fiber, of type $I_6^*$. We will prove that for fixed $\deg(\Li)$ there the two possible configurations can be realized as elliptic surface.

\begin{Prop} Suppose $\pi: X \rightarrow C$ is an elliptic surface with one singular fiber. The fiber is of type $I_{12k-6}^*$ or $I_{12k}$ for some $k\in \Z_{>0}$, and $g(C)\geq k$ in the first case and $g(C) \geq k+1$ in the second. Conversely all these configurations occur.\end{Prop}

\begin{proof} Since the Euler characteristic of the singular fiber is $12\deg(\Li)$, the only possible configurations are $I_\nu$ and $I_\nu^*$.

  Fix a rational elliptic surface $\pi: X \rightarrow \Ps^1$  with 3 fibers of type $III$, and one fiber of type $I_3$. (The existence follows from \cite{Per}.)

  Fix $k$ a positive integer. Take a curve $C$ such that $\phi: C \rightarrow \Ps^1$, has degree $4k-2$, and is ramified at the four points which have the singular fibers, and above such a point there is exactly one point.

  The base change $\pi': X\times_{\Ps^1} C \rightarrow C$ has 3 fibers of type $I_0^*$ and one fiber of type $I_{12k-6}$. Twisting by all four points with a singular fiber yields an elliptic surface with one singular fiber and this fiber is of type $I_{12k-6}^*$.

  If we replace $4k-2$ by $4k$, we obtain an elliptic surface with one fiber of type $I_{12k}$.

For any elliptic surface with only one singular fiber and that fiber is of type $I_0^*$, we have the following: the $j$-map $C \rightarrow \Ps^1$ has degree $12k-6$, one point above $\infty$, at most $4k-2$ points above 0, and at most $6k-3$ points above 1728. This implies that the base curve has genus at least $k$. A similar argument shows that in the case that the only singular fiber is of type $I_{12k}$,  the base curve has genus at least $k+1$.
\end{proof}

\end{document}